\documentclass[english]{smfart}
\usepackage[francais,english]{babel}
\usepackage{smfthm}
\usepackage{basics}
\usepackage{verbatim}

\author{Nick Ramsey}
\address{Department of Mathematics, DePaul University}
\email{naramsey@gmail.edu}
\title[]{$B_{\mathrm{Sen}}$ via distributions on weight space}

\begin{document}
\frontmatter
\maketitle
\tableofcontents

\section{Introduction}
In an unpublished preprint \cite{kisinppmf}, Kisin introduces some
$p$-adic period rings with an eye towards capturing periods for a
certain class of overconvergent $p$-adic modular forms.  The
Fontaine-style functors associated to these rings are intended to be a
sort of Betti realization for these forms.  With this in mind, Kisin
essentially ``deforms'' $B_{HT}$ to force non-integral $p$-adic powers
of the cyclotomic character to turn up.  The rings he defines are
closely related to the Iwasawa algebra and non-canonically contain a
copy of $B_{HT}$.

The $p$-adic modular forms he considers are of a limited sort.  In
particular, they are $p$-adic limits of classical forms.  As such,
their weights are $p$-adic limits of classical weights, which is why
$p$-adic powers of the cyclotomic character suffice for his purposes.
One might naturally be led to ask what can be said outside of this
realm of the eigencurve.  The first thing one might try to do is to
further deform $B_{HT}$ to ``see'' all points of $p$-adic weight space
$\scr{W}$.  Given the interpretation of the Iwasawa algebra in terms
of distributions on $\Z_p$, one natural approach is to consider a
suitable ring of distributions on the rigid-analytic space $\scr{W}$.
In so doing one quickly encounters a number of unpleasant features,
including zero-divisors in the ring of distributions and too many
Galois invariants.  Both of these seem to be related to the presence
of torsion in $\scr{W}$ (which conveniently lies just outside the
region containing the $\Z_p$-powers of the cyclotomic character considered
by Kisin).

This led the author to consider the quotient of $\scr{W}$ by its
torsion.  We show that this quotient has the structure of a rigid
space and that a certain collection of distributions on this space does
indeed furnish a nice ring of periods.  In fact, this ring turns out
to be canonically isomorphic via a ``Fourier transform'' to the the
ring $B_{\mathrm{Sen}}$ introduced by Colmez in \cite{colmezsen}.

\section{Notation}

Fix an odd prime $p$. In this paper, $\scr{W}$ will denote $p$-adic
weight space.  This is a rigid space over $\Q_p$ whose points with
values in a complete extension $L/\Q_p$ are given by $$\scr{W}(L) =
\Hom_{\mathrm{cont}}(\Z_p^\times, L^\times)$$ Let $\tau$ denote the
Teichmuller character $\Z_p^\times\longrightarrow \mu_{p-1}$ and
define $\ip{x} = x/\tau(x)\in 1+p\Z_p$.  Any $\psi\in \scr{W}(K)$ can
be factored as $$\psi(x) = \psi(\tau(x))\psi(\ip{x}) =
\tau(x)^i\psi(\ip{x})$$ for a unique $i\in \Z/(p-1)\Z$.  The space
$\scr{W}$ is the disjoint union of $p-1$ disks $\scr{W}^i$ associated
to these values of $i$.

Let $\gamma$ be a generator of the pro-cyclic subgroup $1+p\Z_p$ (for
example, one could take $\gamma=1+p$).  Then the restriction of $\psi$
to $1+p\Z_p$ is determined by $\psi(\gamma)$.  Moreover, it is not
difficult to show that $\psi(\gamma)=1+t$ where $t\in L$ satisfies
$|t|<1$ and that moreover any such $t$ defines a character of
$1+p\Z_p$ that we will denote by $\psi_t$.  Thus the $L$-valued points
of $\scr{W}$ are in bijective correspondence with $p-1$ copies of the
open unit disk in $L$ via $$(i,t)\longleftrightarrow \tau^i\psi_t$$

For each $n\geq 0$, let $\scr{W}^0_n$ denote the admissible affinoid in
$\scr{W}^0$ defined by the inequality $$|\psi(\gamma)-1|\leq
p^{-1/p^{n-1}(p-1)}$$ Note that this these affinoids do not depend on
the choice of $\gamma$.  The numbering is set up so that if $\zeta$ is
a primitive $p^n$-th root of unity and $\psi$ is satisfies
$\psi(\gamma) = \zeta$, then $\psi\in \scr{W}^0_n\setminus
\scr{W}^0_{n-1}$.  That is, $\scr{W}^0_n\cap \scr{W}^0_{\mathrm{tors}}
=\scr{W}^0[p^n]$.

In what follows $K$ will denote a complete extension of $\Q_p$.  If
$K/\Q_p$ is finite, then for each integer $n\geq 0$ we define $K_n =
K(\mu_{p^n})\subset \overline{\Q_p}$ and $K_\infty = \cup_n
K_n\subseteq \overline{\Q_p}$.  The letter $\chi$ will always denote
the cyclotomic character

We recall some standard facts about the $p$-adic logarithm and
exponential for future use. See \cite{washingtonbook} for a nice
account of this material.  Let $$\log_p(1+x) = \sum_{k\geq 1}
(-1)^{k+1}\frac{x^k}{k}$$ and $$\exp_p(x) = \sum_{k\geq
  0}\frac{x^k}{k!}$$ These series have radii of convergence $1$ and
$p^{-1/(p-1)}$, respectively.  If $|x|<p^{-1/(p-1)}$, then
$|\log_p(1+x)| = |x|$ and we have $$\exp_p(\log_p(1+x)) = x$$
and $$\log_p(\exp_p(x))=x$$ Also, the identity $$\log_p(uv) =
\log_p(u)+\log_p(v)$$ holds whenever it makes sense ($|u-1|,|v-1|<1$).
One consequence of the last property is that $\log_p(\zeta) = 0$ for
any $p$-power root of unity $\zeta$.  In fact, the converse of this
statement holds as well in the sense that the $p$-power roots of unity
are the only roots of $\log_p$ in its domain of convergence.

\section{The quotient $\scr{W}/\scr{W}_{\mathrm{tors}}$}

The space $\scr{W}$ has $(p-1)$-torsion given by the powers $\tau^i$
of the Teichmuller character.  This torsion is all rational over
$\Q_p$ and the quotient of $\scr{W}$ by this torsion subgroup is
canonically identified (via Teichmuller) with the identity component
$\scr{W}^0$.  The space $\scr{W}^0$ has only $p$-power torsion.

We wish to describe the quotient
$\scr{W}^0/\scr{W}^0_{\mathrm{tors}}$.  Looking only at the
$\C_p$-points of these spaces, the properties of $\log_p$ described
above ensure that it defines an injective map
\begin{eqnarray*}
\scr{W}^0(\C_p)/\scr{W}^0_{\mathrm{tors}}(\C_p) & \longrightarrow &\C_p \\
\psi & \longmapsto & \log_p(\psi(\gamma))
\end{eqnarray*}
It is not difficult to see that this map is also surjective (see the
proof of Theorem \ref{quotstructure} below).  This description has two
deficiencies. The first is that it is only a description on points and
does not identify $\scr{W}^0/\scr{W}^0_{\mathrm{tors}}$ as a rigid
space.  The second is that it depends on the choice of $\gamma$.  The
latter is easily remedied as follows.  Identify $\scr{W}^0$ with the
disk defined by $|t-1|<1$ using $\gamma$ as explained in the previous section,
and define an analytic function $\vartheta$ on $\scr{W}^0$
by $$\vartheta = \frac{1}{\log_p(\gamma)}\log_p(1+t)$$ In terms of
characters, this functions is given by $$\vartheta(\psi) =
\frac{\log_p(\psi(\gamma))}{\log_p(\gamma)}$$ It follows immediately
from the latter description and properties of $\log_p$ that the
function $\vartheta$ does not depend on the choice of $\gamma$ and
thus furnishes a canonical analytic function on $\scr{W}^0$ that is
defined over $\Q_p$ (say, by taking $\gamma=1+p$) and
satisfies $$\vartheta(\psi\psi') = \vartheta(\psi)+\vartheta(\psi')\ \ \mbox{for all}\ \ \psi, \psi'\in \scr{W}^0$$

The first difficulty is more subtle, as the notion of quotient in the
setting of rigid spaces is not well-known.  We shall adopt the
following definition of quotient in our context (see
\cite{conradtemkin} for a much more general discussion of such
quotients).  Let $G$ be a rigid-analytic group (such as $\scr{W}^0$ or
$\scr{W}^0_n$) and $H$ be a subgroup (such as
$\scr{W}^0_{\mathrm{tors}}$ or $\scr{W}^0_n[p^n]$, respectively) and
let $H\times G \rightrightarrows G$ denote the equivalence relation
defined by multiplication and projection. We will say that a rigid
space $X$ equipped with a surjective map $G\longrightarrow X$ is the
quotient of $G$ by $H$ if the compositions $$H\times
G\rightrightarrows G\longrightarrow X$$ coincide and the resulting
map $$H\times G \longrightarrow G\times_X G$$ is an isomorphism.

\begin{lemm}\label{thetabound}
  If $\psi\in \scr{W}^0_n$, we have $|\vartheta(\psi)|\leq
  p^{n-1/(p-1)}$.
\end{lemm}
\begin{proof}
If $\psi\in \scr{W}^0_0$ then $|\psi(\gamma)-1|\leq
p^{-p/(p-1)}<p^{-1/(p-1)}$, so 
$$|\vartheta(\psi)| =
\left|\frac{\log_p(\psi(\gamma))}{\log_p(\gamma)}\right|=
p|\log_p(\psi(\gamma))|\leq p^{1-p/(p-1)}=p^{-1/(p-1)}$$ which is the claim for
$n=0$.

  Suppose that the claim holds for some $n$ and let
  $\psi\in\scr{W}_{n+1}^0$.  Observe that $$|\psi(\gamma)^p-1| =
  |((\psi(\gamma)-1)+1)^p-1| = |
  (\psi(\gamma)-1)^p+p(\psi(\gamma)-1) (1+\cdots)|\leq
  p^{-1/p^{n-1}(p-1)}$$ so that $\psi^p\in \scr{W}^0_n$.
  Thus $$p^{-1}|\vartheta(\psi)|=|p\vartheta(\psi)|=|\vartheta(\psi^p)|\leq
  p^{n-1/(p-1)}$$ which establishes the claim for $n+1$ and thus for
  all $n$ by induction.
\end{proof}

The upshot of this lemma is that $\vartheta$ defines a canonical
analytic function 
\begin{equation}\label{thetan}
\vartheta:\scr{W}^0_n\longrightarrow \A^1_{p^{n-1/(p-1)}}
\end{equation}
 for each $n\geq 0$, where $\A^1_{p^{n-1/(p-1)}}$ denotes the affinoid
 ball of radius of $p^{n-1/(p-1)}$ centered at $0$ in the rigid-analytic
 affine line $\A^1$.

\begin{theo}\label{quotstructure}
The map $\vartheta$ in (\ref{thetan}) identifies
$\A^1_{p^{n-1/(p-1)}}$ with the quotient of $\scr{W}^0_n$ by
$\scr{W}^0_n[p^n]$
\end{theo}
\begin{proof}
  Let us first show that the map $\vartheta$ in (\ref{thetan}) is
  surjective.  Suppose that $x\in \C_p$ has $|x|\leq p^{n-1/(p-1)}$.
  Then $$|p^n\log_p(\gamma)x| \leq p^{-n-1}p^{n-1/(p-1)} =
  p^{-1-1/(p-1)}<p^{-1/(p-1)}$$ so $y=\exp_p(p^n\log_p(\gamma)x)$ is
  defined.  Let $z\in \C_p$ satisfy $z^{p^n}=y$ and let $\psi$ be the
  unique point of $\scr{W}^0$ with $\psi(\gamma)=z$.  Then $|z-1|\leq
  p^{-1/p^{n-1}(p-1)}$, as is easy to see using an inductive binomial
  theorem argument akin to the one in the proof of Lemma
  \ref{thetabound}.  Thus $\psi\in \scr{W}^0_n$
  and $$p^n\vartheta(\psi)=
  \frac{\log_p(\psi(\gamma)^{p^n})}{\log_p(\gamma)}
  =\frac{\log_p(\exp_p(p^n\log_p(\gamma)x))}{\log_p(\gamma)} =p^n x$$
  so $\vartheta(\psi)=x$ and $\vartheta$ is surjective.

  Let $A$ denote the affinoid algebra of $\scr{W}_n^0$ which we will
  identify with the Tate algebra of power series in $t$ with
  coefficients in $\Q_p$ that are strictly convergent for $|t|\leq
  p^{-1/p^{n-1}(p-1)}$.  That is, $$A = \left\{\left.\sum
  a_kt^k\ \right|\ |a_k|p^{-k/p^{n-1}(p-1)}\to 0
  \ \mbox{as}\ k\to\infty\right\}$$ Similarly, we let $B$ denote the
  affinoid algebra of $\A^1_{p^n-1/(p-1)}$, which we identify with the
  power series over $\Q_p$ in the variable $s$ that are strictly
  convergent for $|s|\leq p^{n-1/(p-1)}$.  The map $\vartheta$
  corresponds to a function $B\longrightarrow A$, and by suggestive
  abuse of notation, we will denote the image of $s$ under this
  pull-back by $\vartheta(1+t)$.  Explicitly, $$\vartheta(1+t) =
  \frac{\log_p(1+t)}{\log_p(\gamma)} =
  \frac{1}{\log_p(\gamma)}\sum_{k\geq 0}(-1)^{k+1}\frac{t^k}{k}$$
  which lies in $A$, as is evident by taking $\gamma=1+p$, for example.

  Note that the functional properties of $\log_p$ correspond to (in
  fact, result from) formal properties of this series, which we use
  with only minor comments below.

  The two maps $\scr{W}_n^0[p^n]\times\scr{W}_n^0\rightrightarrows
  \scr{W}_n^0$ comprising the equivalence relation correspond to the
  maps
  \begin{eqnarray*}
    A & \longrightarrow & A/((1+t)^{p^n}-1)\ \widehat{\otimes}\ A \\
    t & \longmapsto & 1\otimes t \\
    t & \longmapsto & 1\otimes t + t\otimes 1 + t\otimes t
  \end{eqnarray*}
We claim that the compositions of these maps with the map
$B\longrightarrow A$ coincide.  Indeed
\begin{eqnarray*}
  \vartheta(1+(1\otimes t +t\otimes 1+t\otimes t)) &=&
  \vartheta((1+1\otimes t)(1+t\otimes 1)) \\ &=& \vartheta(1+1\otimes
  t)+\vartheta(1+t\otimes 1) \\ &=& \vartheta(1+1\otimes
  t)+p^{-n}\vartheta((1+t\otimes 1)^{p^n}) \\ &=& \vartheta(1+1\otimes
  t) 
\end{eqnarray*}
where the last equality follows because $$\vartheta((1+t\otimes
1)^{p^n}) = \vartheta(1+ ( (1+t\otimes 1)^{p^n}-1))$$ lies in the
ideal generated by $(1+t\otimes 1)^{p^n}-1$, which is $0$.  The upshot
of the equality of these two functions is that we have a well-defined map
\begin{eqnarray*}
  A\widehat{\otimes}_B A & \longrightarrow &
  A/((1+t)^{p^n}-1)\ \widehat{\otimes}\ A \\
  1\otimes t & \longrightarrow & 1\otimes t \\
  t\otimes 1 & \longrightarrow & 1\otimes t + t\otimes 1 + t\otimes t
\end{eqnarray*}
 To complete the proof, we must show that this map is
an isomorphism.  Let us try to define an inverse map via
\begin{eqnarray*}
  A/((1+t)^{p^n}-1)\ \widehat{\otimes}\ A &\longrightarrow &
  A\widehat{\otimes}_B A\\ 1\otimes t & \longmapsto & 1\otimes t
  \\ t\otimes 1 &\longmapsto & \frac{1+t\otimes 1}{1+1\otimes t}-1
  =(1+t\otimes 1)\left(\sum_{k\geq 0}(-1)^k(1\otimes t)^k\right) - 1
\end{eqnarray*}
If this map is well-defined, then it is a simple formal matter to
check that it is the inverse of the map defined above.  To see that it
is well-defined, first note that the series in the definition has
radius of convergence $1$ and thus lies in $A$.  Finally, we must show
that $(1+t\otimes 1)^{p^n}-1$ gets sent to $0$ under the above proposed
map.  Thus, we need to check $$\left(\frac{1+t\otimes
  1}{1+1\otimes t}\right)^{p^n}-1=0$$ in $A\widehat{\otimes}_BA$,
which is in turn equivalent to $(1+t)^{p^n}\otimes 1 = 1\otimes
(1+t)^{p^n}$ in $A\widehat{\otimes}_BA$.  By definition of the tensor
product, this will follow if we can find an element $b\in B$ whose
image in $A$ under $B\longrightarrow A$ is $(1+t)^{p^n}$.  We claim
that $$b = \exp_p(p^n\log_p(\gamma)s) = \sum_{k\geq
  0}\frac{(p^n\log_p(\gamma))^k}{k!}s^k$$ is such an element.  First,
since $\exp_p(x)$ is strictly convergent on $|x|\leq p^{-1-1/(p-1)}<
p^{-1/(p-1)}$, the given series is strictly convergent for $|s|\leq
p^{n-1/(p-1)}$ and thus lies in $B$.  The image of $b$ under the map
$B\longrightarrow A$ is $$\exp_p(p^n\log_p(\gamma)\vartheta(1+t)) =
\exp_p\left(p^n \log_p(\gamma)
\frac{\log_p(1+t)}{\log_p(\gamma)}\right) = (1+t)^{p^n}$$ as desired.
\end{proof}

\begin{defi}
  For each $n\geq 0$, let $X_n$ denote the quotient
  $\scr{W}_n^0/\scr{W}_n^0[p^n]$.  By the previous theorem, $X_n$ is
  an affinoid defined over $\Q_p$ equipped with a canonical
  isomorphism $\vartheta:X_n\stackrel{\sim}{\longrightarrow}
  \A^1_{p^{n-1/(p-1)}}$.
\end{defi}

Theorem \ref{quotstructure} has the following immediate consequence.

\begin{coro}\label{maincoro}
  The ring $\O(X_n\widehat{\otimes}K)$ consists of functions of the
  form $$f = \sum a_k\vartheta^k$$ with $a_k\in K$ and
  $|a_k|p^{k(n-1/(p-1))}\longmapsto 0$.  In terms of this expansion we
  have $$\|f\|_{\sup{}} = |a_k|p^{k(n-1/(p-1))}$$ for all such $f$.
\end{coro}

For each $n\geq 0$ the natural map
\begin{equation}\label{trans}
X_n\longrightarrow X_{n+1}
\end{equation}
 is injective and identifies $X_n$ with an admissible affinoid open in
 $X_{n+1}$.  Gluing over increasing $n$, we conclude that the function
 $\vartheta$ furnishes a canonical
 isomorphism $$\scr{W}^0/\scr{W}^0_{\mathrm{tors}}\stackrel{\sim}{\longrightarrow}
 \A^1$$

\section{Distributions}

To the affinoids $X_n$, we associate spaces of bounded distributions
as follows.
\begin{defi}
  Let $K/\Q_p$ be a complete extension.  A bounded $K$-valued
  distribution on $X_n$ is a $\Q_p$-linear map
  $$\mu:\O(X_n)\longrightarrow K$$ such that
  there exists $C$ such that $|\mu(f)|\leq C \|f\|_{\sup{}}$ for all $f\in
  \O(X_n)$.  We denote the space of all such
  distributions by $\scr{D}(X_n,K)$.
  The map (\ref{trans}) induces a map 
  \begin{equation}\label{disttrans}
  \scr{D}(X_n,K)\longrightarrow \scr{D}(X_{n+1},K)
  \end{equation}
 for each $n$, and we define $$\scr{D}(X_\infty,K) = \varinjlim_n
 \scr{D}(X_n,K)$$
\end{defi}
The maps in (\ref{disttrans}) are injective for all $n\geq 0$, so the
injective limit is simply a union.  This injectivity is not obvious
from the definition, but is a simple consequence of Lemma
\ref{chardist} below.

\begin{rema}\label{extension}
  If $K/L$ is an extension of complete extensions of $\Q_p$, then any
  $\mu\in \scr{D}(X_n,K)$ induces an $L$-linear
  map $$\O(X_n\widehat{\otimes} L)\longrightarrow K$$ by extension of
  scalars in the obvious manner.
\end{rema}

The following lemma allows us to characterize our distributions on $X_n$
via their ``moments'' under the isomorphism of Theorem \ref{quotstructure}.

\begin{lemm}\label{chardist}\
  \begin{enumerate}
  \item If $\mu\in \scr{D}(X_n,K)$, then
    $|\mu(\vartheta^k)|p^{-k(n-1/(p-1))}$ is bounded in $k$.
    
  \item Conversely, if $x_k\in K$ is a sequence such that
    $|x_k|p^{-k(n-1/(p-1))}$ is bounded, then there exists a unique
    $\mu \in \scr{D}(X_n,K)$ such that $\mu(\vartheta^k)=x_k$.
    \end{enumerate}
\end{lemm}
\begin{proof}\
  \begin{enumerate}
  \item Since $\mu\in \scr{D}(X_n,K)$, there exists $C$ such that
    $|\mu(f)|\leq C\|f\|_{\sup{}}$ for all $f\in \O(X_n)$.
    Thus $$|\mu(\vartheta^k)|\leq C\|\vartheta^k\|_{\sup{}}=
    Cp^{k(n-1/(p-1))}$$

  \item Let $x_k$ be as in the statement. For $f\in
    \O(X_n)$ we may write write $$f = \sum_k
    a_k\vartheta^k$$ with $|a_k|p^{k(n-1/(p-1))}\longrightarrow 0$ by
    Corollary \ref{maincoro}.  The hypotheses on $a_k$ and $x_k$
    ensure that $a_kx_k\longrightarrow 0$, so we may define $\mu(f) =
    \sum_k a_kx_k$.  Note that 
    \begin{eqnarray*}
    |\mu(f)|\leq \sup_k |a_k||x_k| &=&
    \sup_k (|a_k|p^{k(n-1/(p-1))}) (|x_k|p^{-k(n-1/(p-1))}) \\ &\leq&
    \|f\|_{\sup{}}\cdot \sup_k |x_k|p^{-k(n-1/(p-1))}
    \end{eqnarray*}
    again by Corollary \ref{maincoro}.  Thus $\mu$ is bounded and
    defines an element of $\scr{D}(X_n,K)$.
    The uniqueness of $\mu$ follows because $\mu(\sum b_k\vartheta^k)
    = \sum b_k\mu(\vartheta^k)$ holds for any $\mu\in \scr{D}(X_n,K)$
    and any $\sum b_k\vartheta^k\in \O(X_n)$ by the boundedness of
    $\mu$, so any such $\mu$ is \emph{determined} by its moments.
  \end{enumerate}
\end{proof}

The group structure on $X_n$ endows  $\scr{D}(X_n,K)$ with a
convolution product.  To define this product, we will need a lemma.
For $f\in \O(X_n)$ and $\varphi\in X_n$, define a function $T_\varphi
f$ on $X_n$ by $$T_\varphi f(\psi)=f(\varphi\psi)$$ Note that if
$\varphi$ is a $K$-valued point of $X_n$, then $T_\varphi f$ is
naturally an element of $\O(X_n\widehat{\otimes} K)$.
\begin{lemm}
  Let $f\in \O(X_n)$ and let $\mu\in \scr{D}(X_n,K)$.  The
  function $$\varphi\longmapsto \mu(T_\varphi f)$$
is an element of $\O(X_n\widehat{\otimes}K)$
\end{lemm}
\begin{proof}
  First note that this function makes sense by Remark \ref{extension}
  and the comment preceding the lemma.  By Corollary \ref{maincoro},
  $f$ may be written as $f = \sum a_k\vartheta^k$ with
  $|a_k|p^{k(n-1/(p-1))}\longrightarrow 0$.  We have 
    \begin{eqnarray*}
      (T_\varphi f)(\psi) = f(\varphi\psi) &=& \sum_k a_k
      \vartheta(\varphi\psi)^k \\ &=& \sum_k
      a_k(\vartheta(\varphi)+\vartheta(\psi))^k\\ & =& \sum_k
      \sum_{m\leq k}
      a_k\binom{k}{m}\vartheta(\varphi)^{m}\vartheta(\psi)^{k-m} \\
\end{eqnarray*}
    Thus, since $\mu$ is bounded we have
    $$\mu(T_\varphi f) = \sum_k \sum_{m\leq k}a_k
    \binom{k}{m}\mu(\vartheta^{k-m})\vartheta(\varphi)^m =
    \sum_m\left(\sum_{k\geq
      m}a_k\binom{k}{m}\mu(\vartheta^{k-m})\right)\vartheta(\varphi)^m$$
    According to Corollary \ref{maincoro}, in order to complete the proof we
    must show that $$\left|\sum_{k\geq m} a_k\binom{k}{m}
    \mu(\vartheta^{k-m})\right|p^{m(n-1/(p-1))}\longrightarrow 0$$ as
    $m\longrightarrow \infty$.  This expression is bounded by the
    supremum over $k\geq m$
    of $$|a_k||\mu(\vartheta^{k-m})|p^{m(n-1/(p-1))} =
    |a_k|p^{k(n-1/(p-1))} |\mu(\vartheta^{k-m})|p^{(m-k)(n-1/(p-1))}$$
    By Corollary \ref{maincoro}, there exists $C$ such that
    $|\mu(\vartheta^{k-m})|\leq Cp^{(k-m)(n-1/(p-1))}$ for all $k\geq
    m$, so the previous bound is at most $C|a_k|p^{k(n-1/(p-1))}$.
    This quantity tends to zero in $k$, and therefore in $m$ since
    $m\leq k$.
\end{proof}

For $\mu,\nu\in \scr{D}(X_n,K)$ we define the \emph{convolution of
  $\mu$ and $\nu$} by $$(\mu*\nu)(f) =
\mu(\varphi\longmapsto \nu(T_\varphi f))$$ This is well-defined by
the previous lemma, and clearly defines an element of
$\scr{D}(X_n,K)$.  It is trivial to check that the distribution
$\mu_{\bf{1}}\in\scr{D}(X_0,\Q_p)$ given by $\mu_{{\bf{1}}}(f) =
f(\bf{1})$, that is, the Dirac distribution associated to the identity
character, is a two-sided identity for the convolution product.  As we
will see below, the rings $\scr{D}(X_n,K)$ are in fact commutative
integral domains.

\begin{exem}\label{diracconv}
  Let $\psi,\psi'\in \scr{W}^0$.  The convolution of the Dirac
  distributions associated to these two characters
  is $$(\mu_\psi*\mu_{\psi'})(f) = \mu_\psi(\varphi\longmapsto
  \mu_{\psi'}(T_\varphi f)) = \mu_\psi(\varphi\longmapsto
  f(\varphi\psi')) = f(\psi\psi')$$ which is to say
  $\mu_\psi*\mu_{\psi'}=\mu_{\psi\psi'}$.
\end{exem}

\section{Galois action and theta operator}

The group $G=\Gal(\overline{\Q}_p/\Q_p)$ acts on the $\C_p$-valued
points of $X_n$ with $g\in G$ acting on the class of a character
$\psi$ to give the class of $g\circ\psi$.  We define an action of $G$
 on $\O(X_n\widehat{\otimes}\C_p)$ via
$$f^g(\psi) = g(f(g^{-1}\circ \psi))$$ In terms of the expansion of
Corollary \ref{maincoro} this action is simply given by the action of
$G$ on the coefficients $a_k$ since $\vartheta$ is defined over
$\Q_p$, and hence the isomorphism of Theorem \ref{quotstructure}
intertwines this action with the usual one (given by action on
coefficients) on analytic functions on the ball
$\A^1_{p^{n-1/(p-1)}}$.

Let $G_n=G_{\Q_p,n}\subset G$ be the subgroup $G_n =
\Gal(\overline{\Q_p}/\Q_p(\mu_{p^n}))$. Then for $g\in G_n$ we have
$v(\log_p(\chi(g))) = v(\chi(g)-1)\geq n$ and hence
$$p^{k(n-1/(p-1))}\left|\frac{(\log_p(\chi(g)))^k}{k!}\right|\leq 1$$
It follows from Corollary \ref{maincoro} that
$$\exp_p(\vartheta\log_p(\chi(g))) = \sum_k
\frac{(\log_p(\chi(g)))^k}{k!}\vartheta^k$$ is an analytic function on
$X_n$.  This function allows us to define an
action of $G_n$ on $\scr{D}(X_n,\C_p)$ by the
formula
$$(g\cdot \mu)(f) = g( \mu
(f^{g^{-1}}\cdot\exp_p(\vartheta\log_p(\chi(g)))\ ))$$ It is an easy
matter to check that this preserves the boundedness condition and
defines an action of $G_n$.  The collection 
$\{\scr{D}(X_n,\C_p)\}$ together with the
respective actions by $G_n$ is an instance of a
``$G_{\Q_{p,\infty}}$-module'' in the sense of Colmez in \cite{colmezsen}.
\begin{rema}\label{fullgalois}
  In the special case $n=0$, the above recipe actually defines an
  action the entire group $G$ on $\scr{D}(X_0,\C_p)$ if we agree that
  $\log_p(\chi(g))$ is to be interpreted as
  $\log_p(\chi(g)\tau(\chi(g))^{-1})$ (equivalently, if we agree that
  $\log_p$ be extended to $\Z_p^\times$ so that the $(p-1)^{st}$ roots of
  unity be sent to zero).
\end{rema}

Let us now turn to the $\Theta$ operator.  For
$\mu\in\scr{D}(X_n,\C_p)$, define a distribution $\Theta\mu$ by the
formula $(\Theta\mu)(f) = \mu(f\cdot\vartheta)$.  The operator
$\Theta$ so defined preserves the space of bounded distributions and
evidently commutes with the action of $G_n$.

\begin{exem}
  Let us determine explicitly the action of $G_n$ on the Dirac
  distribution $\mu_\psi$ for $\psi$ a point of $\scr{W}^0$.
  \begin{eqnarray*}
    (g\cdot \mu_\psi)(f) &=& g(\mu_{\psi}(f^{g^{-1}}\cdot
    \exp_p(\vartheta\log_p(\chi(g))))) \\ &=&
    g(f^{g^{-1}}(\psi)\exp_p(\vartheta(\psi)\log_p(\chi(g)))) \\ &=&
    f(g\circ\psi)\exp_p(\vartheta(g\circ\psi)\log_p(\chi(g)))
  \end{eqnarray*}
  so $$g\cdot\mu_\psi = \mu_{g\circ\psi}\cdot
  \exp_p(\vartheta(g\circ\psi)\log_p(\chi(g)))$$

In particular, consider the $\Q_p$-valued character
$\psi(x)=x^k\tau(x)^{-k}\in \scr{W}^0_0$.  By Remark \ref{fullgalois},
it makes sense to consider $g\cdot \mu_\psi$ for any $g\in G$, and the
above computation shows that $$g\cdot \mu_\psi =
\mu_\psi\cdot\exp_p(k\log_p(\chi(g)\tau(\chi(g))^{-1}))=
\mu_\psi\cdot\chi(g)^k\tau(\chi(g))^{-k}$$ Thus the $\C_p$-subalgebra
(= $\C_p$-submodule) of $\scr{D}(X_0,\C_p)$ generated by these
characters for all $k\in \Z$ is rather akin to the ring $B_{HT}$, but
with all characters twisted by a Teichmuller power to lie in
$\scr{W}^0$.  Alternatively, one may choose a nonzero $\alpha\in\C_p$
with the property that $g(\alpha)=\tau(\chi(g))\alpha$.  Then we have
(see Example  \ref{diracconv})
\begin{eqnarray*}
g\cdot(\alpha\mu_{x\tau(x)^{-1}})^k &=&
g\cdot(\alpha^k\mu_{x^k\tau(x)^{-k}})\\ &=&
g(\alpha^k)(g\cdot\mu_{x^k\tau(x)^{-k}})\\ &=&
(\tau(\chi(g))^k\alpha^k)
(\mu_{x^k\tau(x)^{-k}}\chi(g)^k\tau(\chi(g))^{-k})\\ &=& \chi(g)^k
(\alpha^k\mu_{x^k\tau(x)^{-k}})
\end{eqnarray*}
It follows that the $\C_p$-subalgebra (= $\C_p$-submodule) of
$\scr{D}(X_0,\C_p)$ generated by the $\alpha^k\mu_{x^k\tau(x)^{-k}}$
for $k\in \Z$ is isomorphic to $B_{HT}$ as a $G$-module, so we have a
non-canonical injection $B_{HT}\hookrightarrow \scr{D}(X_0,\C_p)$
corresponding to the choice of $\alpha$.  This is reminiscent of the
non-canonical injections of $B_{HT}$ into the rings considered in
\cite{kisinppmf} and \cite{colmezsen}.
\end{exem}

\section{Relation to $B_{\mathrm{Sen}}$}

In \cite{colmezsen}, Colmez introduces the ring $B_{\mathrm{Sen}}$
defined as the collection of power series in $\C_p[\![T]\!]$ with
positive radius of convergence.  This ring has an ascending filtration
$\{B_{\mathrm{Sen}}^n\}$ where $B_{\mathrm{Sen}}^n$ consists of power
series of radius of convergence at least $p^{-n}$.  Colmez defines a
``$G_{\Q_p,\infty}$'' action on this filtered ring, meaning a
compatible collection of actions of $G_n$ on $B_{\mathrm{Sen}}^n$, by
acting in the natural way on coefficients and setting $g\cdot T = T +
\log_p(\chi(g))$.

\begin{defi}
Let $\mu\in \scr{D}(X_n,\C_p)$.  
The \emph{Fourier transform} of $\mu$ is the formal series
$$\scr{F}_n(\mu) = \mu(\exp(T\vartheta) ) := \sum_k
\frac{\mu(\vartheta^k)}{k!}T^k$$ 
\end{defi}
Using the usual estimate for the 
$p$-divisibility of $k!$ we have
$$\left|\frac{\mu(\vartheta^k)}{k!}T^k\right|\leq |\mu(\vartheta^k)|
p^{-k(n-1/(p-1))}(p^n|T|)^k$$ If $|T|<p^{-n}$, then this tends to $0$
as $k\to\infty$ by Lemma \ref{chardist}, so the series defining
$\scr{F}_n(\mu)$ has radius of convergence at least $p^{-n}$.  Thus
$\scr{F}_n$ is a $\C_p$-linear
map $$\scr{F}_n:\scr{D}(X_n,\C_p)\longrightarrow B_{\mathrm{Sen}}^n$$
\begin{prop}
  The Fourier transform $\scr{F}_n$ is an injective ring homomorphism.
\end{prop} 
\begin{proof}
That $\scr{F}_n$ is injective follows from Corollary \ref{maincoro}.  Observe that 
\begin{eqnarray*}
  \scr{F}_n(\mu*\nu) &=& (\mu*\nu)(\exp(T\vartheta)) \\ &=&
  \mu(\varphi\mapsto \nu(T_\varphi\exp(T\vartheta))) \\ &=&
  \mu(\varphi\mapsto \nu(\exp(T\vartheta(\varphi\psi)) \\ &=&
  \mu(\varphi\mapsto \nu(\exp(T\vartheta(\varphi))
  \exp(T\vartheta(\psi))\ )) \\ &=&
  \mu(\exp(T\vartheta))\nu(\exp(T\vartheta)) =
  \scr{F}_n(\mu)\scr{F}_n(\nu)
\end{eqnarray*}
Lastly, note that $$\scr{F}_n(\mu_{\mathbf{1}}) =
\mu_{\mathbf{1}}(\exp(T\vartheta)) =
\exp(T\vartheta(\mathbf{1}))=\exp(0)=1$$ Thus $\scr{F}_n$ is indeed a ring
homomorphism.
\end{proof}

It follows immediately from this that the rings $\scr{D}(X_n,K)$ are
commutative integral domains.  In fact, we can say rather more.  Let
$P(T)=\sum a_kT^k$ have positive radius of convergence.  Then there
exists a positive integer $m$ such that $|a_k|p^{-km}\longrightarrow
0$.  The standard estimates on the $p$-divisibility of $k!$ imply that
$$|k!a_k|p^{-k(m+1-1/(p-1))} =
|a_k|p^{-km}|k!|p^{-kp/(p-1)}\longrightarrow 0$$ so by Lemma
\ref{chardist} there exists $\mu\in \scr{D}(X_{m+1},\C_p)$ with
$\mu(\vartheta^k) = k!a_k$. Thus we have $\scr{F}_{m+1}(\mu) =
P(T)$, and, although the individual $\scr{F}_n$ are not isomorphisms
(it is not difficult to check that they are not individually
surjective), the limit $$\scr{F}:\scr{D}(X_\infty,\C_p)
=\varinjlim_n\scr{D}(X_n,\C_p) \stackrel{\sim}{\longrightarrow}
B_{\mathrm{Sen}}$$ is an isomorphism of rings.

We claim that the Fourier transform $\scr{F}_n$ intertwines the action
of $G_n$ that we have defined on bounded distributions with the
action defined by Colmez.  The following calculations are somewhat
formal, but can be made rigorous with power series expansions with no
difficulty.  We have
\begin{eqnarray*}
  \scr{F}_n(g\cdot \mu)(T) &=& (g\cdot\mu)(\exp(T\vartheta)) \\ &=&
  g(\mu(\exp(g^{-1}(T)\vartheta)\exp(\vartheta\log_p(\chi(g))))) \\&=&
  g(\mu( \exp((g^{-1}(T)+\log_p(\chi(g)))\vartheta))) \\ &=&
  g(\scr{F}_n(\mu)(g^{-1}(T)+\log_p(\chi(g)))) 
\end{eqnarray*}
which is equal to $g$ applied to $\scr{F}_n(\mu)$ in the sense of
Colmez, since $\log_p(\chi(g))\in \Q_p$ is Galois-invariant.

In a similar manner, we can determine the way in which $\Theta$
interacts with the Fourier transform.
$$\scr{F}_n(\Theta\mu) = (\Theta\mu)(\exp(T\vartheta)) =
\mu(\vartheta\exp(T\vartheta)) = \frac{d}{dT}\scr{F}_n(\mu)$$ Note
that this is the negative of the operator called $\Theta$ in
\cite{colmezsen}.

Recall that for a $G_{K_\infty}$-module $M=\cup_n M_n$ in the sense of
\cite{colmezsen} we define $M^{G_{K_\infty}} = \cup_n M_n^{G_{K_n}}$.
\begin{prop}\
\begin{enumerate}
\item  $\scr{D}(X_\infty,\C_p)^{G_{\Q_{p,\infty}}} = \Q_{p,\infty}$
\item If $V$ is $d$-dimensional $\C_p$-representation of $G$, then 
  $$(V\otimes_{\C_p} \scr{D}(X_\infty,\C_p))^{G_{\Q_p,\infty}} := \varinjlim
  (V\otimes_{\C_p} \scr{D}(X_n,\C_p))^{G_n}$$ is isomorphic to
  $D_{\mathrm{Sen}}(V)$ as a $\Q_{p,\infty}$ vector space equipped
  with an endomorphism $\Theta$.  In particular, it is $d$-dimensional
  over $\Q_{p,\infty}$
\end{enumerate}
\end{prop}
\begin{proof}
  We have already shown that we have a compatible sequence of
  $G_n$-equivariant injections
  $$\scr{F}_n:\scr{D}(X_n,\C_p)\longrightarrow B^n_{\mathrm{Sen}}$$ that
    are an isomorphism in the injective limit.  It follows easily that
    they induce an isomorphism on invariants, so the result follows
    from Th\'eor\`eme 2 of \cite{colmezsen}.
\end{proof}

\begin{rema}
  This result holds more generally with $\Q_p$ replaced by any finite
  extension $K/\Q_p$, as Th\'eor\`eme 2 of \cite{colmezsen} is proven
  in this generality.
\end{rema}

\end{document}